\documentclass[a4paper,12pt]{amsart}
\usepackage{amssymb}
\usepackage{ifthen}
\usepackage{graphicx}
\nonstopmode \numberwithin{equation}{section}
\newtheorem{thm}[equation]{Theorem}
\newtheorem{cor}[equation]{Corollary}
\newtheorem{lem}[equation]{Lemma}
\newtheorem{prop}[equation]{Proposition}

\newtheorem{conj}{Conjecture}

\theoremstyle{definition}

\newtheorem{prob}[equation]{Problem}

\newenvironment{rem}{%
\bigskip
\noindent \textsl{{\sl Remark. }}}{\bigskip}
\newenvironment{rems}{%
\bigskip
\noindent \textsl{{\sl Remarks. }}}{\bigskip}

\newcounter {own}
\def\theown {\thesection       .\arabic{own}}

\newenvironment{pf}[1][]{%
 \vskip 3mm
 \noindent
 \ifthenelse{\equal{#1}{}}%
  {{\slshape Proof. }}%
  {{\slshape #1.} }%
 }%
{\qed\bigskip}

\newcounter{alphabet}
\newcounter{tmp}

\newcommand{\ID}{{\mathbb D}}

\newcommand{\IC}{{\mathbb C}}

\def\be{\begin{equation}}
\def\ee{\end{equation}}

\newcommand{\bee}{\begin{enumerate}}
\newcommand{\eee}{\end{enumerate}}

\newcommand{\blem}{\begin{lem}}
\newcommand{\elem}{\end{lem}}
\newcommand{\bthm}{\begin{thm}}
\newcommand{\ethm}{\end{thm}}
\newcommand{\bcor}{\begin{cor}}
\newcommand{\ecor}{\end{cor}}
\newcommand{\beg}{\begin{examp}}
\newcommand{\eeg}{\end{examp}}
\newcommand{\begs}{\begin{examples}}
\newcommand{\eegs}{\end{examples}}
\newcommand{\bdefe}{\begin{defin}}
\newcommand{\edefe}{\end{defin}}
\newcommand{\bprob}{\begin{prob}}
\newcommand{\eprob}{\end{prob}}
\newcommand{\bei}{\begin{itemize}}
\newcommand{\eei}{\end{itemize}}

\newcommand{\bcon}{\begin{conj}}
\newcommand{\econ}{\end{conj}}
\newcommand{\bcons}{\begin{conjs}}
\newcommand{\econs}{\end{conjs}}
\newcommand{\bprop}{\begin{prop}}
\newcommand{\eprop}{\end{prop}}
\newcommand{\br}{\begin{rem}}
\newcommand{\er}{\end{rem}}
\newcommand{\brs}{\begin{rems}}
\newcommand{\ers}{\end{rems}}
\newcommand{\bo}{\begin{obser}}
\newcommand{\eo}{\end{obser}}
\newcommand{\bos}{\begin{obsers}}
\newcommand{\eos}{\end{obsers}}
\newcommand{\bpf}{\begin{pf}}
\newcommand{\epf}{\end{pf}}
\newcommand{\ba}{\begin{array}}
\newcommand{\ea}{\end{array}}
\newcommand{\beq}{\begin{eqnarray}}
\newcommand{\beqq}{\begin{eqnarray*}}
\newcommand{\eeq}{\end{eqnarray}}
\newcommand{\eeqq}{\end{eqnarray*}}

\newcommand{\ds}{\displaystyle}

\begin{document}
\bibliographystyle{amsplain}
\title{Region of variability for
spirallike functions with respect to a boundary point}
\author{S. Ponnusamy}
\address{S. Ponnusamy, Department of Mathematics,
Indian Institute of Technology Madras, Chennai-600 036, India.}
\email{samy@iitm.ac.in}
\author{A. Vasudevarao}
\address{A. Vasudevarao, Department of Mathematics,
Indian Institute of Technology Madras, Chennai-600 036, India.}
\email{alluvasu@iitm.ac.in}
\author{M. Vuorinen}
\address{M. Vuorinen, Department of Mathematics,
FIN-20014 University of Turku, Finland.}
\email{vuorinen@utu.fi}
%%{\bf 2000

\date{\today
%March. 06, 06
; File: pvdev$_{-}$6final.tex}
\begin{abstract}

Let ${\mathcal F}_{\mu} $ denote the class of all non-vanishing
analytic functions $f$ in the unit disk $\mathbb{D}$
with $f(0)=1$, and for $\mu\in\IC$,
such that ${\rm Re\,}{\mu}>0$ satisfying
%${\rm Re\,}P_f(z) >0$ in  ${\mathbb D}$, where
$$ {\rm Re\,} \left (\frac{2\pi}{\mu} \frac{zf'(z)}{f(z)}+ \frac{1+z}{1-z}\right )
>0 \quad \mbox{ in  ${\mathbb D}$}
.
$$
For any fixed $z_0$ in the unit disk
%, $a\in\mathbb{C}$ with $|a|\leq 1$ 
and $\lambda\in\overline{\mathbb{D}}$,
we shall determine the region of variability $V(z_0,\lambda)$ for
$\log f(z_0)$ when $f$ ranges over the class
\beqq
\mathcal{F}_{\mu}(\lambda)  & = &
\left\{\frac{}{} f\in{\mathcal F}_{\mu} :\,
f'(0)=\frac{\mu}{\pi}(\lambda-1) 
%~\quad\mbox{ and }~ \right.\\
%&& \quad\left.
%f''(0)=\frac{\mu}{\pi}\left(a(1-|\lambda|^2)+\frac{\mu}{\pi}
%(\lambda-1)^2-(1-{\lambda}^2)\right)
\right\}.
\eeqq
In the final section we graphically illustrate the region of variability
for several sets of parameters.

\end{abstract}

\thanks{The first author was supported by NBHM (DAE, sanction No. 48/2/2006/R\&D-II),
while the second author was supported by
NBHM (DAE, sanction No. 48/2/2006/R\&D-II) and CIMO
(Grant no.15.5.2007/TM-07-5076/CIMO Fellowship),
Academy of Finland, and  Research project ``Quasiconformal Mappings'' of Matti Vuorinen}

\maketitle
\pagestyle{myheadings}
\markboth{S. Ponnusamy, A. Vasudevarao and M. Vuorinen}{Regions of variability}

\section{Introduction}

We denote the class of analytic functions in the unit disk
$\ID=\{z\in\IC:\,|z|<1\}$ by ${\mathcal H}(\ID)$, and we think of
$\mathcal{H}(\ID)$ as a topological vector space endowed with the
topology of uniform convergence over compact subsets of
$\mathbb{D}$. Denote by $ {\mathcal S}^*$ the subclass of functions
$\phi\in{\mathcal H}(\ID)$ with $\phi(0)=0$ such that $\phi$ maps
$\mathbb{D}$ univalently onto a domain $\Omega=\phi(\mathbb{D})$
that is starlike with respect to the origin. That is,
$t\phi(z)\in\phi(\mathbb{D})$ for each  $t\in [0,1]$.
It is well known that for
$\phi\in\mathcal{H}(\mathbb{D})$ with $\phi(0)=0=\phi'(0)-1$, $\phi
\in {\mathcal S}^*$ if and only if
$${\rm Re\,}\left (\frac{z\phi '(z)}{\phi (z)}\right )>0, \quad z\in \ID.
$$
Functions in ${\mathcal S}^*$ are referred to as starlike
functions. Denote by $ {\mathcal C}$ the subclass of functions
$\phi\in{\mathcal H}(\ID)$ with $\phi(0)=0$ such that $\phi$ maps
$\mathbb{D}$ univalently  onto a convex domain. It is well known
that for $\phi\in\mathcal{H}(\mathbb{D})$ with
$\phi(0)=0=\phi'(0)-1$, $\phi \in {\mathcal C}$ if and only if
$z\phi'\in{\mathcal S}^*$. Functions in ${\mathcal C}$ are
referred to as convex functions. We refer to the books
\cite{Du,Goodman's_book} for a detailed discussion on these two
classes. Although, the class of starlike functions (with respect
to an interior point) has been studied extensively among many
other subclasses, little was known about starlike functions with
respect to a boundary point until the work of Robertson
\cite{Robertson}. Motivated by the work in \cite{Robertson} and
characterizations of this class of functions, some advancement in
this direction has taken place (see
\cite{Silverman-Silvia,lyzzak-84,elin-Reich-01,lecko-01}). On the
other hand, there does not seem to be any development on
spiral-like functions with respect to a boundary point until the
recent work of Elin  et al. \cite{elin-Reich-00} (see also
\cite{elin-Reich-01}). More recently, Aharonov et al.
\cite{Aharonov} provide with a natural geometric approach to
discuss spiral-like functions with respect to a boundary point and
the conditions described in \cite{Aharonov} cover the results
studied by others. On the other hand,  several authors have
studied region of variability problems for various subclasses of
univalent functions in  ${\mathcal H}(\ID)$, see
\cite{Paatero,Pinchuk,
%samy-vasudev-Yanahi1,samy-vasudev-Yanahi4,
samy-vasudev2,samy-vasudev3,samy-vasudev-yan3,Yanagihara1,Yanagihara2}. For example,
it is well-known that for each fixed $z_0\in \ID$, the region of
variability
$$V(z_0)=\{\log \phi ' (z_0):\,\phi \in {\mathcal C},~ \phi'(0)=1 \}
$$
is the set $\{\log (1-z)^{-2}:\,|z|\leq |z_0|\}$.

Let ${\mathcal F}_\mu$ denote the class of functions
$f\in {\mathcal H}(\ID)$, and non-vanishing in $\mathbb{D}$ with $f(0)=1$,
and for $\mu\in\IC$, such that ${\rm Re\,}{\mu}>0$
satisfying
$${\rm Re\,}P_f(z) >0, \quad z\in{\mathbb D},
$$
where \be\label{pvdev-6a-eq1} P_f(z)= \frac{2\pi}{\mu}\frac{zf'(z)}{f(z)}+
\frac{1+z}{1-z}. \ee Clearly $P_f(0)=1$. Basic properties and a number of equivalent
characterizations of the class $\mathcal {F}_\mu $
are formulated in \cite{Aharonov}. The case $\mu =\pi$ coincides with the
class introduced by Robertson \cite{Robertson} who has generated interest
on this class, and its associated classes. It is also known that functions in
$\mathcal{F}_\pi$ are either close-to-convex or just the constant $1$.

% This class was introduced by M.S.
%Robertson\cite{Robertson} and a renewed  intrest has been
%generated on this class $\mathcal{F}$ and .... generalization in
%\cite{}.

For $f\in{\mathcal F}_\mu $, we denote by $\log f$ the single-valued
branch of the logarithm of $f$ with $\log f(0)=0$. Herglotz
representation for analytic function with positive real part in
$\mathbb{D}$ shows that if $f\in\mathcal{F}_\mu $, then there exists a
unique positive unit measure $\nu$ on $(-\pi,\pi]$ such that
$$\frac{2\pi}{\mu}\frac{zf'(z)}{f(z)}+ \frac{1+z}{1-z}=
\int_{-\pi}^{\pi}\frac{1+ze^{-it}}{1-ze^{-it}}\,d\nu(t),
$$
and hence, a computation gives that
$$\log f(z)=\frac{\mu}{\pi}
\int_{-\pi}^{\pi}\log\left(\frac{1-z}{1-ze^{-it}}\right)d\nu(t);
$$
or equivalently
$$f(z)= (1-z)^{\mu/\pi}\exp\left\{\frac{\mu}{\pi}
\int_{-\pi}^{\pi}\log\left(\frac{1}{1-ze^{-it}}\right)d\nu(t)\right\}.
$$
Let ${\mathcal B}_0$ be the class of analytic functions $\omega$
in $\mathbb{D}$ such that $| \omega(z)|\leq 1$ in $\mathbb{D}$ and
$\omega(0)=0$. Consequently, for each $f\in {\mathcal F}_\mu$ there
exists an $\omega_f \in {\mathcal B}_0$ of the form
% we have ${\rm Re\,}P_f(z)>0$ and
%$P_f(0)=\alpha$, the function
\be\label{pvdev-6a-eq2}
\omega_f(z)=\frac{P_f(z)-1}{P_f(z)+1},\quad z\in\mathbb{D},
\ee
and conversely. It is a simple exercise to  see that
\be\label{pvdev-6a-eq2a}
P'_f(0)=2\omega'_f(0)=2\left(\frac{\pi}{\mu}f'(0)+1\right).
\ee
Suppose that $f\in {\mathcal F}_\mu$. Then, a simple application of the
classical Schwarz lemma
(see for example \cite{Du,samy-book3,samy-herb}) shows that
$$|P'_f(0)|= 2|(\pi/\mu)f'(0)+1|\leq 2,
$$
because $|\omega'_f(0)|\leq 1$. Using (\ref{pvdev-6a-eq2}), one can obtain by a
computation that
$$\frac{\omega''_f(0)}{2}=\frac{P_f''(0)}{4}-{\lambda}^2,\quad\mbox{ and }
P_f''(0)=\frac{4\pi}{\mu}f''(0)-\frac{4\mu}{\pi}(\lambda-1)^2+4
$$
so that
$$\frac{\omega''_f(0)}{2}=\frac{\pi}{\mu}f''(0)-\frac{\mu}{\pi}(\lambda-1)^2+1-{\lambda}^2.
$$
Also if we let
$$
g(z)=\frac{\frac{\omega_f(z)}{z}-\lambda}
{1-\overline{\lambda}\frac{\omega_f(z)}{z}},\quad\mbox{ for } |\lambda|<1,
$$
and $g(z)=0$ for $|\lambda|=1$, then we see that
$$g'(0)=\left\{\ba{ll}
\ds \left.\frac{1}{1-|\lambda|^2}\left(\frac{\omega_f(z)}{z}\right)'\right |_{z=0}
=\frac{1}{1-|\lambda|^2}\left(\frac{\omega''_f(0)}{2}\right)
& \mbox{for $|\lambda|< 1$} \\ [6mm]
\ds 0 & \mbox{for $|\lambda|= 1$}.
 \ea \right .
$$
We note that for $|\lambda|<1$,
\beqq
|g'(0)|\leq 1 & \Leftrightarrow & \frac{|\omega''_f(0)|}{2(1-|\lambda|^2)}\leq 1\\
&\Leftrightarrow & \frac{1}{1-|\lambda|^2}\left|\frac{\pi}{\mu}f''(0)
-\frac{\mu}{\pi}(\lambda-1)^2+1-{\lambda}^2 \right|\leq 1\\
& \Leftrightarrow & f''(0)=\frac{\mu}{\pi}\left(a(1-|\lambda|^2)+\frac{\mu}{\pi}
(\lambda-1)^2-(1-{\lambda}^2) \right)
\eeqq
for some $a\in\mathbb{C}$ with $|a|\leq 1$.
Consequently, for
$\lambda\in\overline{\mathbb{D}}=\{z\in\mathbb{C}:\,|z|\leq 1\}$
%, $a\in\mathbb{C}$ with $|a|\leq 1$ 
and for $z_0\in\mathbb{D}$ fixed,
it is natural to introduce
\beqq
\mathcal{F}_{\mu}(\lambda)  & = &
\left\{\frac{}{} f\in{\mathcal F}_{\mu} :\,
f'(0)=\frac{\mu}{\pi}(\lambda-1) 
%~\quad\mbox{ and }~ \right.\\
%&& \quad\left.
%f''(0)=\frac{\mu}{\pi}\left(a(1-|\lambda|^2)+\frac{\mu}{\pi}
%(\lambda-1)^2-(1-{\lambda}^2)\right)
\right\}\\
V (z_0,\lambda) &  = &  \{\log f(z_0):\, f\in {\mathcal F}_\mu(\lambda)\}
\eeqq

From (\ref{pvdev-6a-eq2a}) and the normalization condition introduced
in the class ${\mathcal F}_\mu(\lambda)$, we observe that
$\omega'_f(0)=\lambda$. The main aim of this paper is to determine
the region of variability $V (z_0,\lambda)$ for $\log f(z_0)$ when
$f$ ranges over the class ${\mathcal F}_\mu (\lambda)$. The precise
geometric description of the set $V (z_0,\lambda)$ is established
in Theorem \ref{pvdev-6a-th1}.
%determine explicitly the
%region of variability $V_{\phi}(z_0,\lambda)$ of $ f(z_0)$ when $f$
%ranges over the class $C(\lambda)$.

\section{Basic properties of $V (z_0,\lambda)$ and the  Main result}

To state our main theorem, we need some  preparation. For a
positive integer $p$, let
$$({\mathcal S}^*)^p=\{f=f_0^p:\, f_0\in {\mathcal S}^* \}
$$
and recall the following result from \cite{Yanagihara1}.

\blem\label{pvdev-6a-lem1} Let $f$ be an analytic function in
${\mathbb D}$ with $f(z) = z^p + \cdots $. If
$$ {\rm Re} \,  \left( 1+ z \frac{f''(z)}{f'(z)} \right)> 0 , \quad z \in {\mathbb D} ,
$$
then $f \in ({\mathcal S}^*)^p$. \elem

%Although we could not find any historical reference for a proof of
%Lemma \ref{pvdev-6a-lem1}, it might be well known
%(see  \cite{Goodman's_book,Hallenbeck-Livingston}), and
%We refer to \cite{Yanagihara1} for an analytic proof of Lemma \ref{pvdev-6a-lem1}.
Now, we list down some basic properties of $V (z_0,\lambda)$.
\bprop\label{pvdev-6a-pro01}
We have
\bee
\item $V (z_0,\lambda)$ is compact.
\item $V (z_0,\lambda)$ is convex.
\item for $|\lambda|=1$ or $z_0=0$,
\be\label{pvdev-6a-eq20}
V (z_0,\lambda)=\left \{\frac{\mu}{\pi}\log\left(\frac{1-z_0}{1-\lambda z_0}\right)\right \}.
\ee
\item for $|\lambda|<1$ and $z_0\in\mathbb{D}\setminus\{0\}$, $V (z_0,\lambda)$ has
$(\mu/\pi)\log\left(\frac{1-z_0}{1-\lambda z_0}\right)$ as an interior point.
\eee
\eprop
\bpf
(1) Since ${\mathcal F}_\mu (\lambda)$ is a
compact subset of $\mathcal{H}(\ID)$, it follows that
$V (z_0,\lambda)$ is also compact.

(2) If $f_0,f_1\in {\mathcal F}_\mu (\lambda)$ and $0\leq t \leq 1$, then
the function
$$f_t(z)=\exp\left\{(1-t)\log f_0(z)+t\log f_1(z)\right\}
$$
is evidently in ${\mathcal F}_\mu (\lambda)$. Also, because of the
representation of $f_t$, we see easily that  the set $V (z_0,\lambda)$ is
convex.

(3) If $z_0=0$, (\ref{pvdev-6a-eq20}) trivially holds. If $|\lambda|=|\omega'_f(0)|=1$,
then it follows from the classical Schwarz lemma that $\omega_f(z)=\lambda z$, which
implies
$$P_f(z)=\frac{1+\lambda z}{1-\lambda z}
~\mbox{ and }~f(z)=\left(\frac{1-z}{1-\lambda z}\right)^{\frac{\mu}{\pi}}.
$$
Consequently,
$$V (z_0,\lambda)
=\left \{\frac{\mu}{\pi}\log\left(\frac{1-z_0}{1-\lambda z_0}\right)\right \}.
$$

(4) For $|\lambda|<1$, and $a\in\overline{\mathbb{D}}$, we define
$$\delta(z,\lambda) = \frac{z+\lambda}{1+\overline{\lambda}z},
$$
and
\be\label{pvdev-6a-eq4}
H_{a,\lambda}(z)  =
\exp\left(\frac{\mu}{\pi}\int_0^z\frac{\delta(a\zeta, \lambda)-1}
{(1-\delta(a\zeta, \lambda)\zeta)(1-\zeta)}\,d\zeta \right),\quad
z\in\mathbb{D}.
\ee
First we claim that $H_{a,\lambda}\in {\mathcal F}_\mu (\lambda)$.
For this, we compute
\begin{eqnarray*}
\frac{2\pi}{\mu}\frac{zH'_{a,\lambda}(z)}{H_{a,\lambda}(z)}
& = & \frac{2z(\delta(az,\lambda)-1)}{(1-\delta(az,\lambda)z)(1-z)}\\
& = & \frac{2z\delta(az,\lambda)}{1-\delta(az,\lambda)z}-\frac{2z}{1-z}
\end{eqnarray*}
and so, we see easily that
$$\frac{2\pi}{\mu}\frac{zH'_{a,\lambda}(z)}{H_{a,\lambda}(z)}+\frac{1+z}{1-z}
=\frac{1+\delta(az,\lambda)z}{1-\delta(az,\lambda)z}.
$$
As $\delta(az,\lambda)$ lies in the unit disk $\mathbb{D}$,
$H_{a,\lambda}\in {\mathcal F}_\mu (\lambda)$ and the claim follows.
Also we observe that
\be\label{pvdev-6a-eq5}
\omega_{H_{a,\lambda}}(z)=z\delta(az,\lambda).
\ee

Next we claim that the mapping ${\mathbb D} \ni a\mapsto \log H_{a,\lambda}(z_0) $
is a non-constant analytic function of $a$
for each fixed $z_0 \in {\mathbb D} \backslash \{ 0 \}$ and
$\lambda\in\mathbb{D}$. To do this, we put
\begin{eqnarray*}
h(z) & = & \left .\frac{2\pi}{\mu(1-|\lambda|^2)}
\frac{\partial}{\partial a}\left\{\frac{}{}\log H_{a,\lambda} (z)\right\}\right |_{a=0}.
%\\
\end{eqnarray*}
A computation gives
\begin{eqnarray*}
h(z) & = &  2\int_0^z \frac{\zeta}{(1-\lambda\zeta)^2}\,d\zeta=z^2+\cdots
\end{eqnarray*}
from which it is easy to see that
$${\rm Re} \,\left\{\frac{zh''(z)}{h'(z)}\right\}=
{\rm Re} \,\left\{\frac{1+\lambda z}{1-\lambda z}\right\}> 0,
\quad z\in\mathbb{D}.
$$
By Lemma \ref{pvdev-6a-lem1} there exists a function $h_0\in
{\mathcal S}^*$  with $h=h_0^2$. The univalence of $h_0$ together with
the condition $h_0(0)=0$ implies that $h(z_0)\neq 0$ for $z_0 \in {\mathbb
D}\setminus \{0\}$. Consequently, the mapping ${\mathbb D} \ni
a\mapsto \log H_{a,\lambda}(z_0) $ is a non-constant analytic
function of $a$ and hence, it is an open mapping. Thus,
$V (z_0,\lambda)$ contains the open set $\{\log
H_{a,\lambda}(z_0):\, |a|<1\}$. In particular,
$$\log H_{0,\lambda}(z_0)=(\mu/\pi)\log \left(\frac{1-z_0}{1-\lambda z_0}\right)
$$
is an interior point of $\{\log H_{a,\lambda}(z_0):\,a\in\mathbb{D}\}\subset V (z_0,\lambda)$.
\epf

We remark that, since $V (z_0,\lambda)$ is a compact convex subset of
$\mathbb{C}$ and has nonempty interior, the boundary
$\partial{V (z_0,\lambda)}$ is a Jordan curve and $V (z_0,\lambda)$
is the union of $\partial{V (z_0,\lambda)}$ and its inner domain.

Now we state our main result and the proof will be presented in
Section \ref{sec2}.

\bthm \label{pvdev-6a-th1} For $\lambda\in\mathbb{D}$ and
$z_0\in\mathbb{D}\setminus \{0\}$, the boundary
$\partial{V (z_0,\lambda)}$ is the Jordan curve given by
\begin{eqnarray*}
%\label{pvdev-6a-eq7}
(-\pi,\pi]\ni \theta \mapsto \log H_{e^{i\theta},\lambda}(z_0) & =
& \frac{\mu}{\pi}\int_0^{z_0}\frac{\delta(e^{i\theta}\zeta, \lambda)-1}
{(1-\delta(e^{i\theta}\zeta, \lambda)\zeta)(1-\zeta)} \,d\zeta.
\end{eqnarray*}
If  $\log f(z_0)= \log H_{e^{i\theta},\lambda}(z_0)$ for some
$f\in{\mathcal F}_\mu (\lambda)$ and $\theta\in (-\pi,\pi]$, then
$f(z)=H_{e^{i\theta},\lambda}(z)$. \ethm

\section{Proof of Theorem \ref{pvdev-6a-th1}}\label{sec2}

\bprop\label{pvdev-6a-pro1} For $f\in {\mathcal F}_\mu (\lambda)$ we
have \be\label{pvdev-6a-eq8}
\left|\frac{f'(z)}{f(z)}-\frac{\mu}{\pi}c(z,\lambda)\right|\leq \frac{|\mu|}{\pi}
r(z,\lambda),
\quad z\in\mathbb{D},
\ee
where
\begin{eqnarray*}
c(z,\lambda)& = & \frac{|z|^2(\overline{z}-\lambda)(1-\overline{\lambda})-
(1-\lambda)(1-\overline{\lambda}\overline{z})}
{(1-z)(1-|z|^2)(1+|z|^2-2{\rm Re\,}({\lambda} z))}, ~~ \mbox{ and }\\
r(z,\lambda) & = & \frac{(1-|\lambda|^2)|z|}
{(1-|z|^2)(1+|z|^2-2{\rm Re\,}({\lambda} z))}.
\end{eqnarray*}
For each $z\in\mathbb{D}\setminus\{0\}$, equality holds if and
only if $f=H_{e^{i\theta},\lambda}$ for some
$\theta\in\mathbb{R}$.
\eprop \bpf Let $f\in{\mathcal F}_\mu (\lambda)$. Then there exists
$\omega_f \in {\mathcal B}_0$ satisfying
(\ref{pvdev-6a-eq2}). As noticed in the introduction through (\ref{pvdev-6a-eq2a}) and
the normalization of $f$, we have $\omega'_f(0)=\lambda$. It follows from
the Schwarz lemma (see for example
\cite{Du,%pommerenke,
samy-book3,samy-herb}) that
\be\label{pvdev-6a-eq9}
\left|\frac{\frac{\omega_f(z)}{z}-\lambda}
{1-\overline{\lambda}\frac{\omega_f(z)}{z}}\right|\leq|z|, \quad z\in\mathbb{D}.
\ee
From  (\ref{pvdev-6a-eq1}) and (\ref{pvdev-6a-eq2}) this is
equivalent to
\be\label{pvdev-6a-eq10}
\left|\frac{\frac{f'(z)}{f(z)}-\frac{\mu}{\pi}A(z,\lambda)}
{\frac{f'(z)}{f(z)}+\frac{\mu}{\pi}B(z,\lambda)}
\right| \leq |z|\, |\tau(z,\lambda)|,
\ee
where
\be\label{pvdev-6a-eq11}
\left\{
\ba{lll}
A(z,\lambda)& =& \ds \frac{\lambda-1}{(1-\lambda z)(1-z)}\\
B(z,\lambda) & = & \ds \frac{1-\overline{\lambda}}{(1- z)(z-\overline{\lambda})}\\
\tau(z,\lambda)& = & \ds \frac{z-\overline{\lambda}}{1-\lambda z}.
\ea
\right.
\ee
A simple calculation shows that the inequality
(\ref{pvdev-6a-eq10}) is equivalent to
\be\label{pvdev-6a-eq12}
\left|\frac{f'(z)}{f(z)}-\frac{\mu}{\pi}\frac{A(z,\lambda)+|z|^2\,
|\tau(z,\lambda)|^2 B(z,\lambda)}
{1-|z|^2\,|\tau(z,\lambda)|^2}\right| \leq \frac{|\mu|}{\pi}
\frac{|z|\,|\tau(z,\lambda)|\,|A(z,\lambda)+B(z,\lambda)|\,}{1-|z|^2\,|\tau(z,\lambda)|^2}.
\ee
Using (\ref{pvdev-6a-eq11}) we can easily see that
\begin{eqnarray}
\label{pvdev-6a-eq13} 1-|z|^2\, |\tau(z,\lambda)|^2 & = &
\frac{(1-|z|^2)(1+|z|^2-2{\rm Re\,}(\lambda z))}{|1-\lambda
z|^2}\nonumber,
\\
\label{pvdev-6a-eq14}
A(z,\lambda)+ B(z,\lambda)
 & = & \frac{1-|\lambda|^2}{(1-\lambda z)(z-\overline{\lambda})}\nonumber
\end{eqnarray}
and
\begin{eqnarray*}
%\hspace{0.75cm}
A(z,\lambda)+|z|^2|\tau(z,\lambda)|^2 B(z,\lambda)
& = & \frac{(\lambda-1)(1-\overline{\lambda}\overline{z})
+|z|^2(\overline{z}-\lambda)(1-\overline{\lambda})}
{(1-z)|1-\lambda z|^2}\nonumber.
\end{eqnarray*}
Thus, by a simple computation, we see that
$$\frac{A(z,\lambda)+|z|^2|\tau(z,\lambda)|^2 B(z,\lambda)}
{1-|z|^2|\tau(z,\lambda)|^2} = c(z,\lambda)
$$
and
$$\frac{|z|\,|\tau(z,\lambda)|\,|A(z,\lambda)+B(z,\lambda)|}{1-|z|^2|\tau(z,\lambda)|^2}
=r(z,\lambda).
$$
Now the inequality (\ref {pvdev-6a-eq8}) follows from these
equalities and (\ref{pvdev-6a-eq12}).

It is easy to see that the equality occurs for a $z\in\mathbb{D}$ in
(\ref {pvdev-6a-eq8}), when $f=H_{e^{i\theta},\lambda}$ for some
$\theta\in\mathbb{R}$. Conversely if the equality occurs for some
$z\in\mathbb{D}\setminus\{0\}$ in (\ref {pvdev-6a-eq8}), then the
equality must hold in (\ref{pvdev-6a-eq9}). Thus from the Schwarz
lemma there exists a $\theta\in\mathbb{R}$ such that
$\omega_f(z)=z\delta(e^{i\theta}z,\lambda)$ for all $z\in\mathbb{D}$. This
implies $f=H_{e^{i\theta},\lambda}$. \epf

The choice of $\lambda=0$ gives the following result which may need
a special mention.

\bcor\label{pvdev-6a-cor01}
For $f\in{\mathcal F}_\mu (0)$ we have
%\be\label{pvdev-6a-eq22}
$$\left|\frac{f'(z)}{f(z)}-\frac{\mu(|z|^2\overline{z}-1)}{\pi(1-z)(1-|z|^4)}\right|\leq
\frac{|\mu|\,|z|}{\pi(1-|z|^4)}, \quad z\in\mathbb{D}.
$$
%\ee

For each
$z\in\mathbb{D}\setminus\{0\}$, equality holds if and only if
$f=H_{e^{i\theta},0}$ for some $\theta\in\mathbb{R}$.
\ecor

%If $f\in{\mathcal F}_\mu (0)$, then by (\ref{pvdev-6a-eq22}) one has
%$$(1-|z|^4)\left|\frac{f'(z)}{f(z)}\right|<
%\frac{|\mu|}{\pi}\left(\frac{1+|z|}{1-|z|}-|z|^2\right).
%$$

\bcor\label{pvdev-6a-cor1}
Let $\gamma:\,z(t)$, $0\leq t\leq 1$, be
a $C^1$-curve in $\mathbb{D}$ with $z(0)=0$ and $z(1)=z_0$. Then
we have
$$V (z_0,\lambda)\subset
%\overline{\mathbb{D}}\left (\frac{\mu}{\pi}C(\lambda, \gamma),
%\frac{|\mu|}{\pi}R(\lambda, \gamma)\right ) =
\left\{w\in\mathbb{C}:\,\left |w-\frac{\mu}{\pi}
C(\lambda, \gamma) \right |\leq \frac{|\mu|}{\pi}R(\lambda, \gamma)\right\},
$$
where
$$C(\lambda, \gamma)=\int_0^1 c(z(t),\lambda)z'(t)\,dt ~\mbox{ and }~ R(\lambda,
\gamma)=\int_0^1 r(z(t),\lambda)|z'(t)|\,dt.
$$
\ecor \bpf
For $f\in {\mathcal F}_\mu(\lambda)$, it follows from
Proposition \ref{pvdev-6a-pro1} that
\begin{eqnarray*}
\left |\log f(z_0)-\frac{\mu}{\pi}C(\lambda, \gamma)\right | & = &
\left|\int_0^1\left\{\frac{f'(z(t))}{f(z(t))}
-\frac{\mu}{\pi}c(z(t),\lambda)\right\}z'(t)\,dt\right|\\
& \leq &
\int_0^1\left|\frac{f'(z(t))}{f(z(t))}-\frac{\mu}{\pi}c(z(t),\lambda)\right|\,
|z'(t)|\,dt \\
& \leq &  \frac{|\mu|}{\pi}\int_0^1 r(z(t),\lambda)|z'(t)|\,dt
= \frac{|\mu|}{\pi}R(\lambda, \gamma).
\end{eqnarray*}
Since $ \log f(z_0)\in V (z_0, \lambda)$ was arbitrary, the
conclusion follows. \epf

For the proof of our next result, we need the following lemma.

\blem\label{pvdev-6a-lem01}
For $\theta\in\mathbb{R}$ and $\lambda\in\mathbb{D}$, the function
$$G(z)=\frac{\mu}{\pi}\int_0^z \frac{ e^{i\theta}\zeta }
{\{1+(\overline{\lambda}e^{i\theta}-\lambda)\zeta-e^{i\theta}{\zeta}^2\}^2}\,
d\zeta, \quad z\in\mathbb{D},
$$
has a double zero at the origin and no zeros elsewhere in
$\mathbb{D}$. Furthermore there exists a starlike univalent
function $G_0$ in $\mathbb{D}$ such that
$G=(\mu/(2\pi))e^{i\theta}G^2_0$ and $G_0(0)= G'_0(0)-1=0$.
\elem \bpf
Let $b = {\rm Im}(\overline{\lambda}e^{i\theta/2}) \in {\mathbb
R}$. Then a computation shows that
\beqq
1+(\overline{\lambda}e^{i\theta}-\lambda)z-e^{i\theta}z^2
%1-\lambda z -e^{i\theta}z
%(z-\overline{\lambda}) &=& 1+(\overline{\lambda} e^{i\theta/2}-\lambda e^{-i\theta /2})
%e^{i\theta /2}z-(e^{i\theta /2}z)^2\\
%&=&1+2ibe^{i\theta /2}z-(e^{i\theta /2}z)^2\\
&=& (1- z/z_1)(1-z/z_2 ), \eeqq
where
$$z_1 =e^{-i\theta /2}(ib+\sqrt{1-b^2}) ~\mbox{ and } ~
z_2=e^{-i\theta /2}(ib-\sqrt{1-b^2}).
$$
From this we have
\beqq
\frac{G''(z)}{G'(z)}-\frac{1}{z}&=& \frac{d}{dz}\left\{\log \frac{G'(z)}{z}\right\}\\
%&=&\frac{d}{dz}\left\{i\theta +\log z - 2\log (1-z/z_1)-2\log (1-z/z_2)\right\}\\
&=&\frac{2/z_1}{1-z/z_1}+\frac{2/z_2}{1-z/z_2} .
\eeqq
Since $|z_1|= |z_2|=1$, we have for $z \in \ID$
$$
{\rm Re} \left\{1+\frac{z G''(z)}{G'(z)}\right\} ={\rm Re}
\left\{\frac{1+z/z_1}{1-z/z_1}\right\} +{\rm Re}
\left\{\frac{1+z/z_2}{1-z/z_2}\right\} > 0 .
$$
Applying Lemma \ref{pvdev-6a-lem1} to $(2\pi/\mu)e^{-i \theta}G(z)$ with
$p=2$ there exists a $G_0 \in S^*$ such that $G=(\mu/(2\pi))e^{i\theta}G^2_0$. \epf

\bprop
Let $z_0\in\mathbb{D}\setminus \{0\}$. Then for
$\theta\in(-\pi,\pi]$ we have
$\log H_{e^{i\theta},\lambda}(z_0)\in\partial V (z_0,\lambda)$.
Furthermore if $\log f(z_0)= \log H_{e^{i\theta},\lambda}(z_0)$
for some $f\in {\mathcal F}_\mu (\lambda)$ and $\theta\in(-\pi,\pi]$,
then $f= H_{e^{i\theta},\lambda}$.
\label{pvdev-6a-pro2}
\eprop\bpf
From (\ref{pvdev-6a-eq4}) we have
\begin{eqnarray*}
\frac{H'_{a, \lambda}(z)}{H_{a, \lambda}(z)}
 & = & \frac{\mu}{\pi}\frac{\delta(az, \lambda)-1}{(1-\delta(az, \lambda)z)(1-z)}\\
 & = & \frac{\mu}{\pi}
 \frac{(\lambda-1)+(1-\overline{\lambda})az}{(1-z)(1+(\overline{\lambda}a-\lambda)z-az^2)}.
\end{eqnarray*}
Using (\ref{pvdev-6a-eq11}) we compute
$$\frac{H'_{a, \lambda}(z)}{H_{a, \lambda}(z)}-\frac{\mu}{\pi}A(z,\lambda)
=\frac{\mu(1-|\lambda|^2)az}{\pi(1-\lambda z)(1+(\overline{\lambda}a-\lambda)z-az^2)},
$$
$$\frac{H'_{a, \lambda}(z)}{H_{a, \lambda}(z)}+\frac{\mu}{\pi}B(z,\lambda) =
\frac{\mu(1-{\lambda}^2)}{\pi(z-\lambda)(1+(\overline{\lambda}a-\lambda)z-az^2)}
$$
and hence we obtain that
%\vspace{8pt}
%$F'_{a, \lambda}(z) -  c(z,\lambda)$
\begin{eqnarray*}
\frac{H'_{a, \lambda}(z)}{H_{a, \lambda}(z)} -  \frac{\mu}{\pi}c(z,\lambda) & = &
\frac{H'_{a, \lambda}(z)}{H_{a, \lambda}(z)}-
\frac{\mu}{\pi}\frac{A(z,\lambda)+|z|^2|\tau(z,\lambda)|^2 B(z,\lambda)}
{1-|z|^2|\tau(z,\lambda)|^2}
\\
 &  =  &
\frac{1}{1-|z|^2 |\tau(z,\lambda)|^2}
\left\{\left(\frac{H'_{a,\lambda}(z)}{H_{a,\lambda}(z)}
-\frac{\mu}{\pi}A(z,\lambda)\right)\right.\\
& & \qquad \left . \frac{}{} -|z|^2|\tau(z,\lambda)|^2
\left(\frac{H'_{a,\lambda}(z)}{H_{a,\lambda}(z)}+\frac{\mu}{\pi}B(z,\lambda)\right)\right\}
\\
& = & \frac{\mu(1-|\lambda|^2)z[a(1-\overline{\lambda}\overline{z})-
\overline{z}(\overline{z}-\lambda)]} {\pi(1-|z|^2)(1+|z|^2-2{\rm Re\,}(\lambda z))
(1+(\overline{\lambda}a-\lambda)z-az^2)}\\
& = & r(z,\lambda)\frac{\mu}{\pi}\frac{az}{|z|} \left
(\frac{|1+(\overline{\lambda}a-\lambda)z-az^2|^2}
{(1+(\overline{\lambda}a-\lambda)z-az^2)^2}\right ).
\end{eqnarray*}
Now by substituting $a=e^{i\theta}$ we easily see that
\begin{eqnarray*}
\frac{H'_{e^{i\theta}, \lambda}(z)}{H_{e^{i\theta}, \lambda}(z)} -
\frac{\mu}{\pi}c(z,\lambda)
& = & r(z,\lambda)\frac{\mu}{\pi}\frac{{e^{i\theta}}z}{|z|} \left
(\frac{|1+(\overline{\lambda}e^{i\theta}-\lambda)z-e^{i\theta}z^2|^2}
{(1+(\overline{\lambda}e^{i\theta}-\lambda)z-e^{i\theta}
z^2)^2}\right ).
\end{eqnarray*}
Putting $G(z)$ as in Lemma \ref{pvdev-6a-lem01}, we get that
\be\label{pvdev-6a-eq15}
\frac{H'_{e^{i\theta},
\lambda}(z)}{H_{e^{i\theta}, \lambda}(z)}- \frac{\mu}{\pi}c(z,\lambda) =
\frac{|\mu|}{\pi}r(z,\lambda)\frac{G'(z)}{|G'(z)|}
\ee
and  there exists a starlike univalent
function $G_0$ in $\mathbb{D}$ such that
$G=(\mu /(2\pi))e^{i\theta}G^2_0$ and $G_0(0)= G'_0(0)-1=0$.
As the function $G_0$ is
starlike, for any $z_0\in\mathbb{D}\setminus\{0\}$ the linear
segment joining $0$ and $G_0(z_0)$ entirely lies in
$G_0(\mathbb{D})$. Now, we define $\gamma_0$ by
\be\label{pvdev-6a-eq16}
\gamma_0:\,z(t)=G_0^{-1}(tG_0(z_0)),
\quad 0\leq t \leq 1.
\ee
Since $G(z(t))=(\mu/(2\pi))e^{i\theta}(G_0(z(t)))^2=
(\mu/(2\pi))e^{i\theta}(tG_0(z_0))^2=t^2G(z_0)$, we have
\be\label{pvdev-6a-eq17}
G'(z(t))z'(t)=2tG(z_0),\quad t\in [0,1].
\ee
Using  (\ref{pvdev-6a-eq17}) and (\ref{pvdev-6a-eq15}) we have
\begin{eqnarray}
\label{pvdev-6a-eq18} \log H_{e^{i\theta},
\lambda}(z_0)-\frac{\mu}{\pi}C(\lambda,\gamma_0)
 & = & \int_0^1\left\{\frac{H'_{e^{i\theta}, \lambda}(z(t))}{H_{e^{i\theta}, \lambda}(z(t))}
-\frac{\mu}{\pi}c(z(t),\lambda)\right\}z'(t)\,dt\\
 & = &\frac{|\mu|}{\pi} \int_0^1 r(z(t),\lambda)
\frac{G'(z(t))z'(t)}{|G'(z(t))z'(t)|}|z'(t)|\,dt\nonumber \\
& = & \frac{G(z_0)}{|G(z_0)|}\frac{|\mu|}{\pi}\int_0^1 r(z(t),\lambda)|z'(t)|\,dt
\nonumber \\
& = & \frac{G(z_0)}{|G(z_0)|}\frac{|\mu|}{\pi}R(\lambda, \gamma_0)\nonumber ,
\end{eqnarray}
where $C(\lambda,\gamma_0)$ and $R(\lambda,\gamma_0)$ are defined as in
Corollary \ref{pvdev-6a-cor1}. Thus, we have
$$\log H_{e^{i\theta},\lambda}(z_0)\in\partial{\overline{\mathbb{D}}}
\left (\frac{\mu}{\pi}C(\lambda, \gamma_0),\frac{|\mu|}{\pi}R(\lambda, \gamma_0)\right ).
$$
Also, from Corollary \ref{pvdev-6a-cor1}, we have
$$\log H_{e^{i\theta},\lambda}(z_0)\in V (z_0,\lambda)\subset\overline{\mathbb{D}}
\left (\frac{\mu}{\pi}C(\lambda,\gamma_0),\frac{|\mu|}{\pi}R(\lambda, \gamma_0)\right ).
$$
Hence, we  conclude that $\log H_{e^{i\theta},\lambda}(z_0)\in
\partial V (z_0,\lambda)$.

Finally, we prove the uniqueness of the curve. Suppose that
$$\log f(z_0)=\log H_{e^{i\theta},\lambda}(z_0)
$$
for some $f\in {\mathcal F}_\mu (\lambda)$ and $\theta\in (-\pi, \pi]$. We introduce
$$h(t)=\frac{\overline{G(z_0)}}{|G(z_0)|}
\left\{\frac{f'(z(t))}{f(z(t))}-\frac{\mu}{\pi}c(z(t),\lambda)\right\}z'(t),
$$
where $\gamma_0:\,z(t)$, $ 0\leq t \leq 1$, is given by
(\ref{pvdev-6a-eq16}). Then, $h(t)$ is continuous function in
$[0,1]$ and satisfies
$$|h(t)|\leq \frac{|\mu|}{\pi}r(z(t),\lambda)|z'(t)|.
$$
Furthermore, we have from (\ref{pvdev-6a-eq18})
\begin{eqnarray*}
%\label{pvdev-6a-eq19}
\int_0^1 {\rm Re} \,h(t)\,dt
 & = &
 \int_0^1 {\rm Re} \, \left\{\frac{\overline{G(z_0)}}{|G(z_0)|}
\left\{\frac{f'(z(t))}{f(z(t))}-\frac{\mu}{\pi}c(z(t),\lambda)\right\}z'(t)\right\}dt
\\
 & = &{\rm Re} \,\left\{\frac{\overline{G(z_0)}}{|G(z_0)|} \left\{\log
f(z_0)-\frac{\mu}{\pi}C(\lambda, \gamma_0)\right\}\right\}\nonumber\\
& = & {\rm Re} \,\left\{\frac{\overline{G(z_0)}}{|G(z_0)|}
\left\{\log
H_{e^{i\theta},\lambda}(z_0)-\frac{\mu}{\pi}C(\lambda,\gamma_0)\right\}\right\}\nonumber
\\
& = & \frac{|\mu|}{\pi}\int_0^1 r(z(t),\lambda)|z'(t)|\,dt\nonumber.
\end{eqnarray*}
Thus, we have
$$h(t)= \frac{|\mu|}{\pi}r(z(t),\lambda)|z'(t)| ~\mbox{ for all $t\in [0,1]$.}
$$
From (\ref{pvdev-6a-eq15}) and (\ref{pvdev-6a-eq17}), it follows
that
$$\frac{f'}{f}=\frac{H'_{e^{i\theta}, \lambda}}{H_{e^{i\theta}, \lambda}}
~\mbox{ on $\gamma_0$.}
$$
By applying the identity theorem for analytic
functions, we get
$$\frac{f'}{f}=\frac{H'_{e^{i\theta}, \lambda}}{H_{e^{i\theta}, \lambda}}
~\mbox{ in $\mathbb{D}$}
$$
%$f'/f=H'_{e^{i\theta}, \lambda}/H_{e^{i\theta}, \lambda}$ in $\mathbb{D}$
and hence, by normalization, $f=H_{e^{i\theta},\lambda}$ in $\mathbb{D}$.
\epf

\noindent {\bf Proof of Theorem \ref{pvdev-6a-th1}.}  We need to prove
that the closed curve
$$(-\pi,\pi]\ni \theta \mapsto \log H_{e^{i\theta},\lambda}(z_0)
$$
is simple. Suppose that
$$\log H_{e^{i\theta_1},\lambda}(z_0)=\log H_{e^{i\theta_2},\lambda}(z_0)
$$
for some $\theta_1,\theta_2\in(-\pi,\pi]$ with $\theta_1\neq\theta_2$.
Then, from Proposition \ref{pvdev-6a-pro2}, we have
$$H_{e^{i\theta_1},\lambda}= H_{e^{i\theta_2},\lambda}.
$$
From (\ref{pvdev-6a-eq5}) this gives a contradiction that
$$e^{i\theta_1}z=\tau\left(\frac{\omega_{H_{e^{i\theta_1},\lambda}}}{z},\lambda\right)
=\tau\left(\frac{\omega_{H_{e^{i\theta_2},\lambda}}}{z},\lambda\right)=e^{i\theta_2}z.
$$
Thus, the curve must be simple.

Since $V (z_0,\lambda)$ is a compact convex subset of $\mathbb{C}$
and has nonempty interior, the boundary $\partial V (z_0,\lambda)$
is a simple closed curve. From Proposition \ref{pvdev-6a-pro1}, the
curve $\partial V (z_0,\lambda)$ contains the curve $(-\pi,\pi]\ni
\theta\mapsto \log H_{e^{i\theta},\lambda}(z_0)$. Recall the fact
that a simple closed curve cannot contain any simple closed curve
other than itself. Thus, $\partial V (z_0,\lambda)$ is given by
$(-\pi,\pi]\ni\theta \mapsto \log H_{e^{i\theta},\lambda}(z_0)$.

\section{Geometric view of Theorem \ref{pvdev-6a-th1}}

  Using Mathematica 4.1,  we describe the boundary of the set
$V(z_0, \lambda)$.  Here we give the Mathematica program
which is used to plot the boundary of the set  $V(z_0, \lambda)$.
We refer \cite{Ruskeepaa} for Mathematica program. The short notations in this program
are of the form: ``z0 for $z_0$'',  ``lam for $\lambda$'' and ``mu for $\mu$''.

%(* The values ``z0, lam, mu'' are for the left hand side picture
%of FIGURE 1 *)
%z0 =-0.173777 + 0.0869191i ;
%lam =-0.196029 + 0.480913i ;
%mu = 32796 + 64560.2i;

\vspace{0.1cm}
{\tt
\begin{verbatim}
Remove["Global`*"];

z0 = Random[]Exp[I* Random[Real, {-Pi, Pi}]]
lam = Random[]Exp[I *Random[Real, {-Pi, Pi}]]
mu = Random[Real, {0, 10^3}] + I *Random[Real, {-10^3, 10^3}]

Q[lam_, the_] := ((lam - 1) + (1 - Conjugate[lam])Exp[I*the]z)/
                 ((1 - z)((1 + ( Conjugate[lam]*Exp[I*the] - lam)*z )
                 - Exp[I*the]*z*z));

myf2[lam_, the_, z0_] := mu/Pi NIntegrate[Q[lam, the], {z, 0, z0}];
image = ParametricPlot[{Re[myf2[lam, the, z0]],
        Im[myf2[lam, the, z0]]}, {the, -Pi, Pi},
        AspectRatio -> Automatic];
(*Clear[z0, lam, mu];*)
\end{verbatim}
}

\vspace{0.5cm}

  The following pictures give the geometric view of the boundary of the
set $V(z_0, \lambda)$. Each of the following figures contain two pictures which describe
the boundary of the set $V(z_0, \lambda)$ for fixed value of $z_0\in\mathbb{D}\setminus\{0\}$,
$\lambda\in\mathbb{D}$  and $\mu\in\mathbb{C}$ such that ${\rm Re \,}\mu>0$. The corresponding values
for each picture are given in a column at the bottom of the picture.
Note that according to  Proposition \ref{pvdev-6a-pro01}
the region bounded by the curve $\partial V(z_0, \lambda)$ is compact and convex.

\begin{figure}[htp]
\begin{center}
\includegraphics[width=5cm]{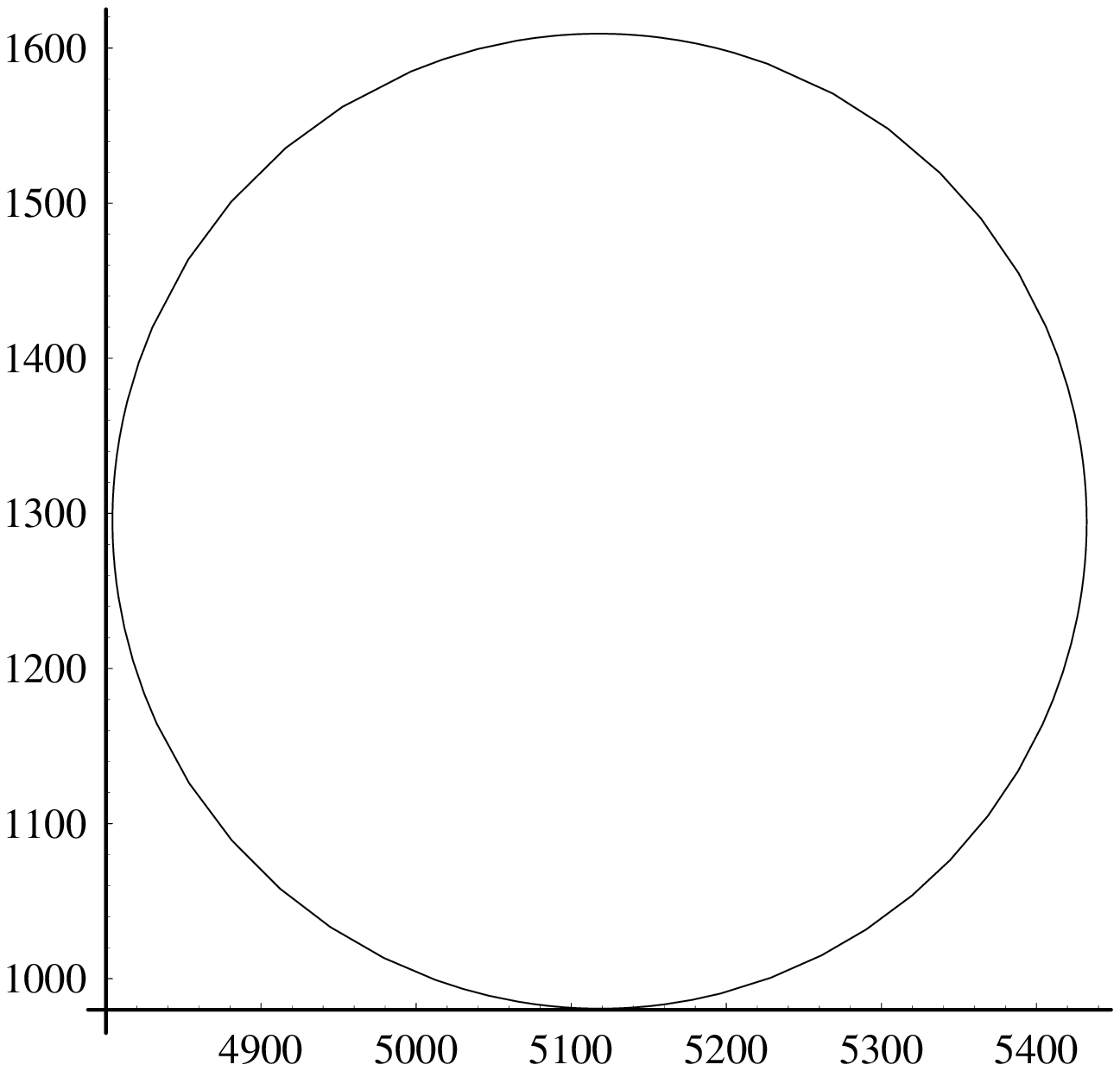}
\hspace{1cm}
\includegraphics[width=5cm]{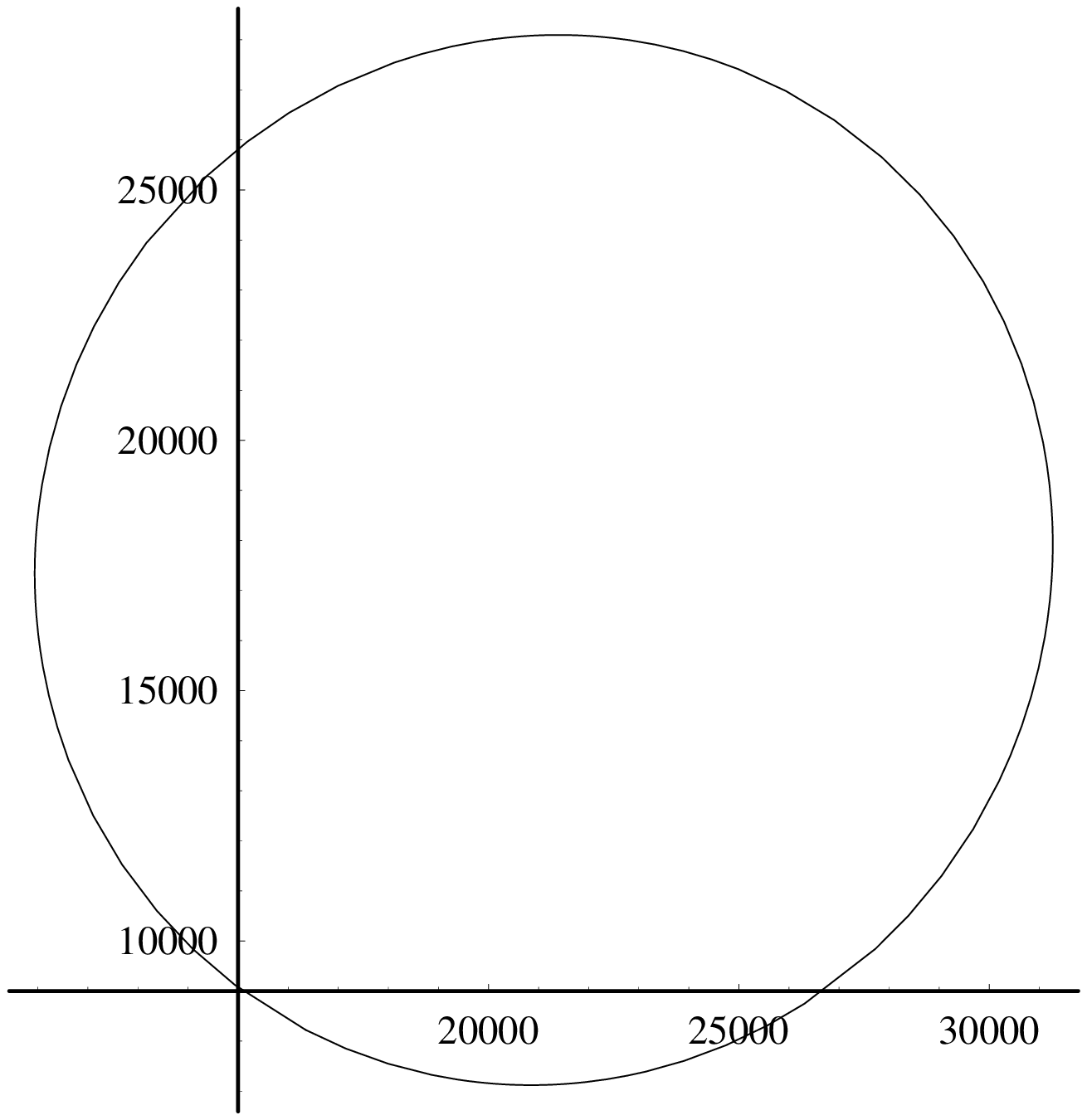}
\end{center}
\hspace{1.5cm}$\partial V(z_0,\lambda)$ \hspace{4.5cm}  $\partial V(z_0,\lambda)$
\caption{Region of variability for $\log f(z_0)$}
\end{figure}
\begin{center}
$\begin{array}{ll}
z_0 =-0.173777 + 0.0869191i       \hspace{3cm} & z_0 =-0.713811 - 0.0997298i    \\
\lambda = -0.196029 + 0.480913i                & \lambda = -0.225338 + 0.323073i  \\
\mu = 32796 + 64560.2i                         & \mu=  69097.4 + 83886.6 i
\end{array}$
\end{center}

\begin{figure}[htp]
\begin{center}
\includegraphics[width=6cm]{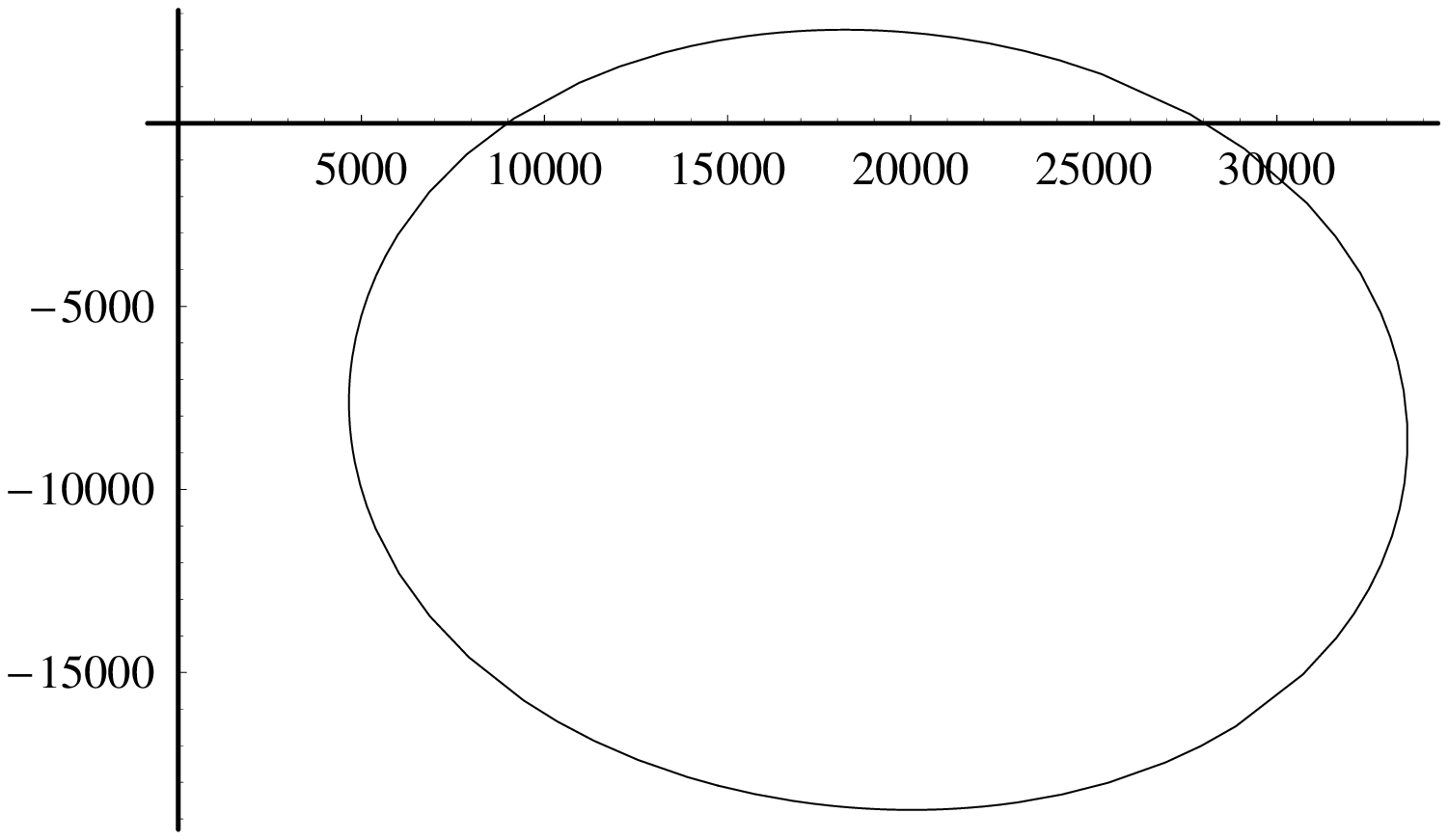}
\hspace{1cm}
\includegraphics[width=5cm]{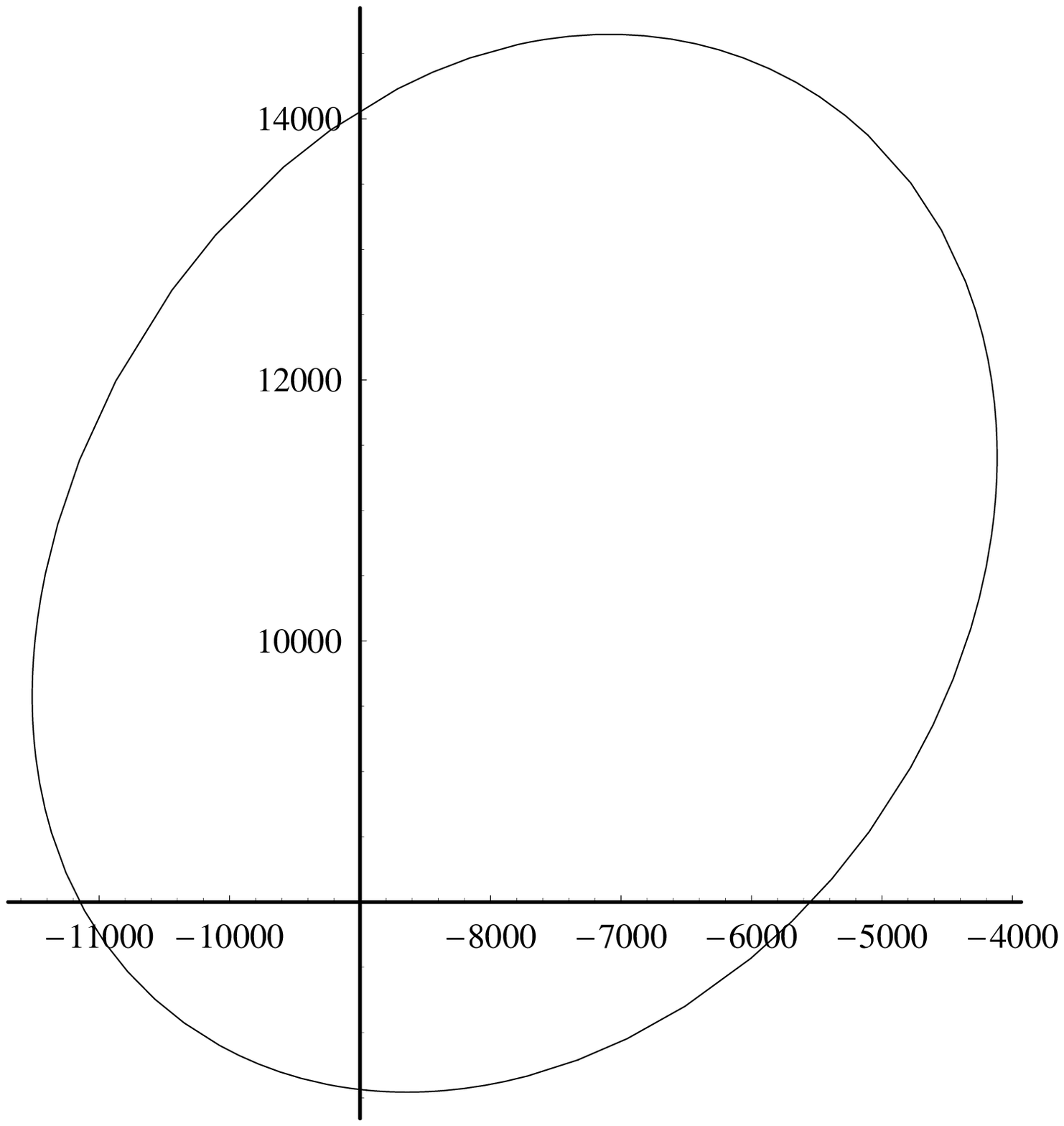}
\end{center}
\hspace{1.5cm}$\partial V(z_0,\lambda)$ \hspace{4.5cm}  $\partial V(z_0,\lambda)$
\caption{Region of variability for $\log f(z_0)$}
\end{figure}
\begin{center}
$\begin{array}{ll}
z_0 =-0.734426 + 0.61942i   \hspace{3cm} & z_0 = -0.69693 - 0.601351i   \\
\lambda = -0.0564481 - 0.00656122i       & \lambda =-0.0416728 - 0.683999i  \\
\mu = 54025 - 5108.28i                   & \mu= 23944.2+ 50613.5i
\end{array}$
\end{center}
\newpage

\begin{figure}[htp]
\begin{center}
\includegraphics[width=5cm]{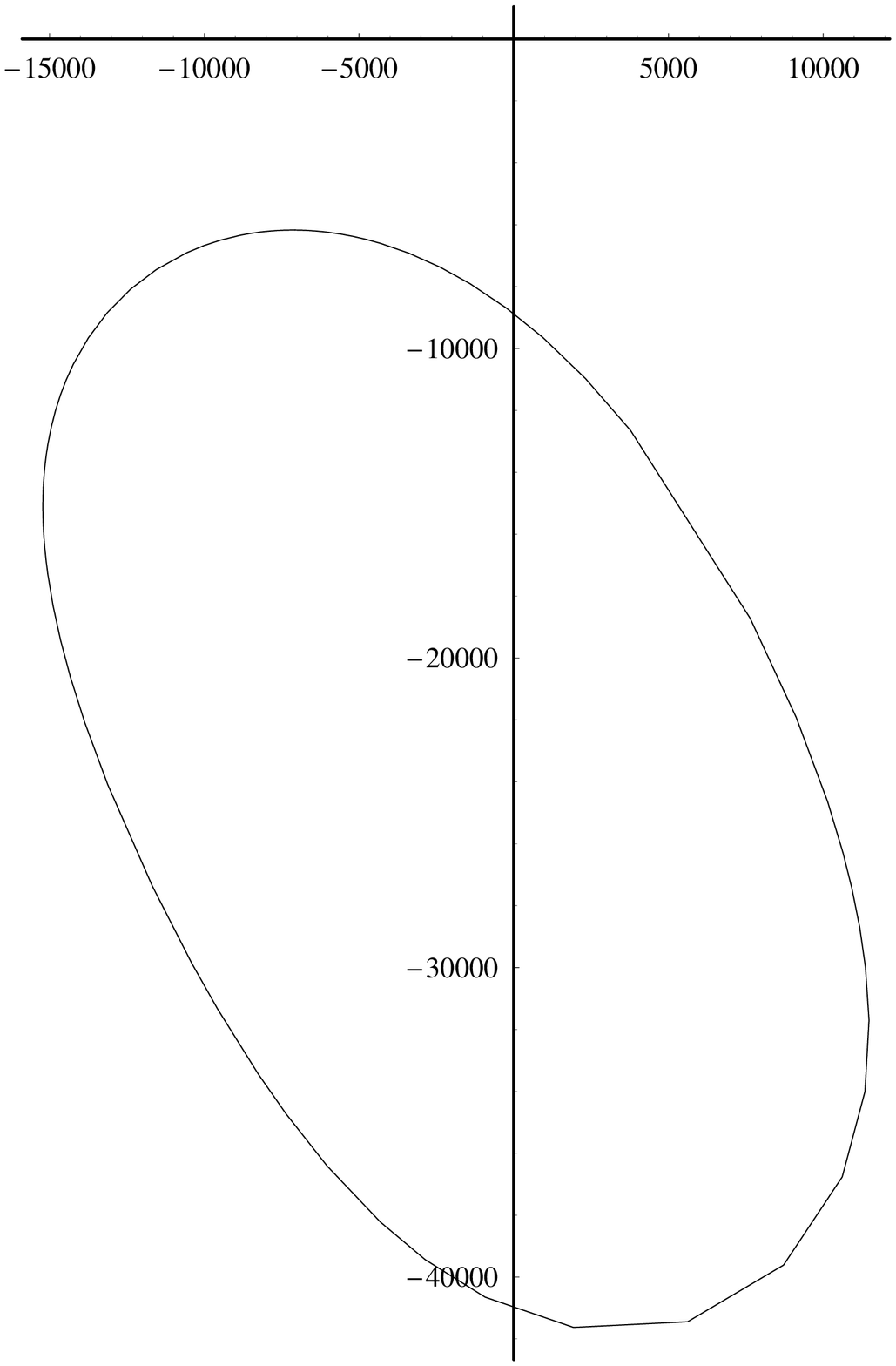}
\hspace{1cm}
\includegraphics[width=5cm]{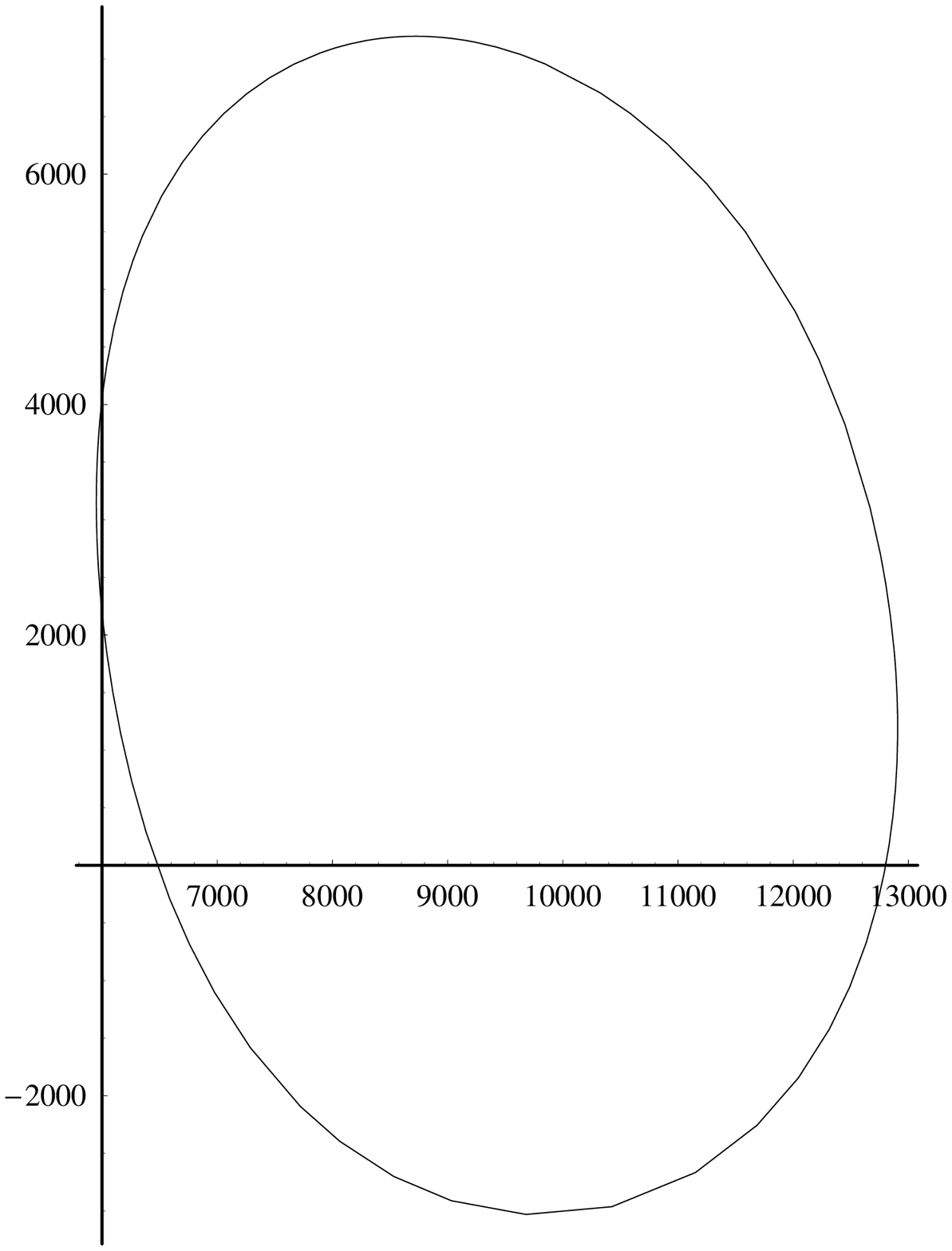}
\end{center}
\hspace{1.5cm}$\partial V(z_0,\lambda)$ \hspace{4.5cm}  $\partial V(z_0,\lambda)$
\caption{Region of variability for $\log f(z_0)$}
\end{figure}
\begin{center}
$\begin{array}{ll}
z_0 = 0.0150249+ 0.994594i      \hspace{3cm} & z_0 = 0.378332 - 0.90135i   \\
\lambda =-0.219752 - 0.256693i               & \lambda =0.366791- 0.600223i  \\
\mu =16828.1- 35690.8i                       & \mu=5006.59- 46769.8i
\end{array}$
\end{center}

\begin{figure}[htp]
\begin{center}
\includegraphics[width=5cm]{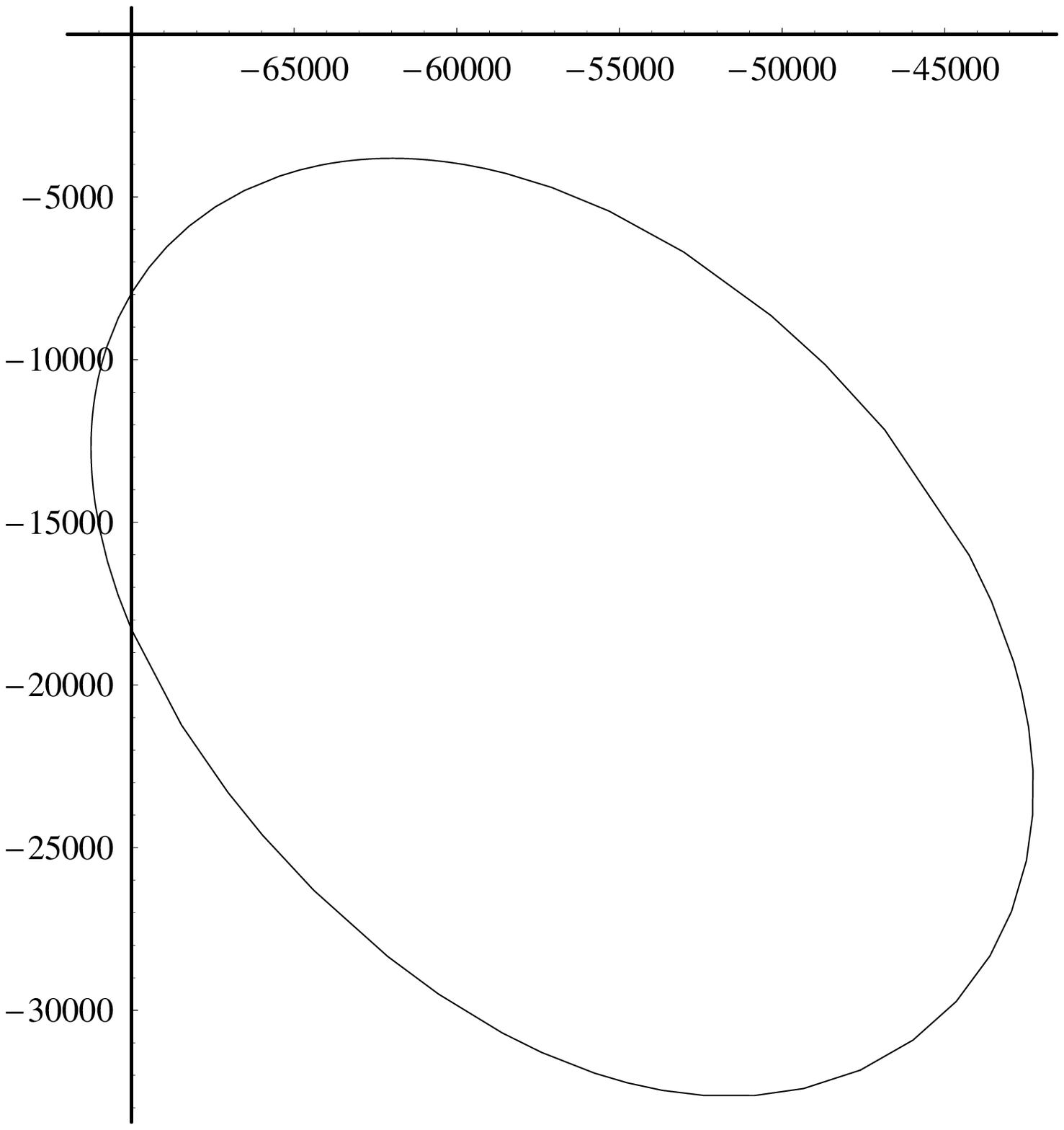}
\hspace{1cm}
\includegraphics[width=6cm]{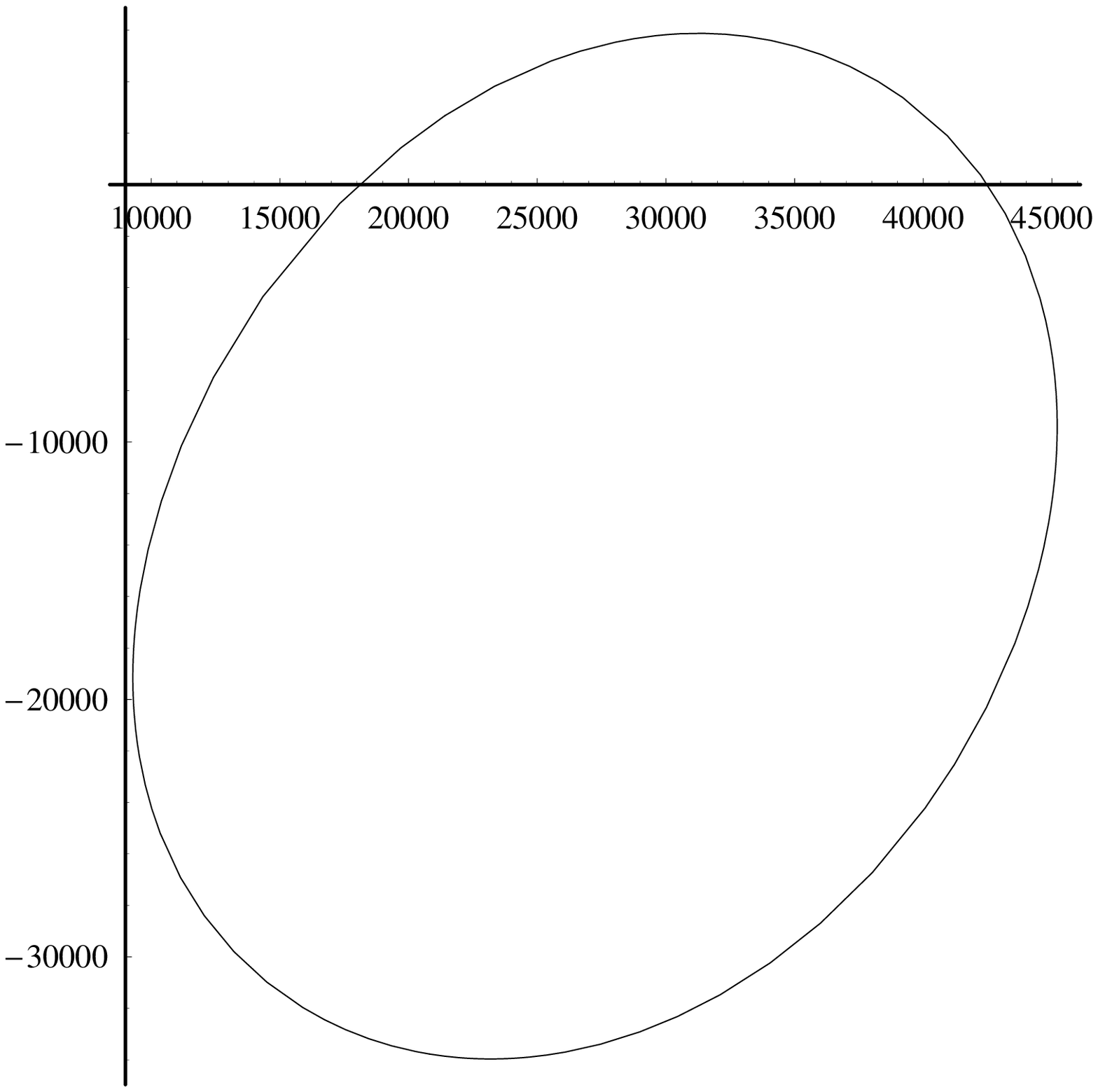}
\end{center}
\hspace{1.5cm}$\partial V(z_0,\lambda)$ \hspace{4.5cm}  $\partial V(z_0,\lambda)$
\caption{Region of variability for $\log f(z_0)$}
\end{figure}
\begin{center}
$\begin{array}{ll}
z_0 =0.80351+ 0.549035i       \hspace{3cm} & z_0 = 0.691568 + 0.644823i   \\
\lambda = -0.55886 + 0.0419296i            & \lambda = 0.126172+ 0.137643i \\
\mu =83278.8 - 90464.3i                    & \mu= 47178.4 + 83497.8i
\end{array}$
\end{center}

\newpage

\begin{figure}[htp]
\begin{center}
\includegraphics[width=5cm]{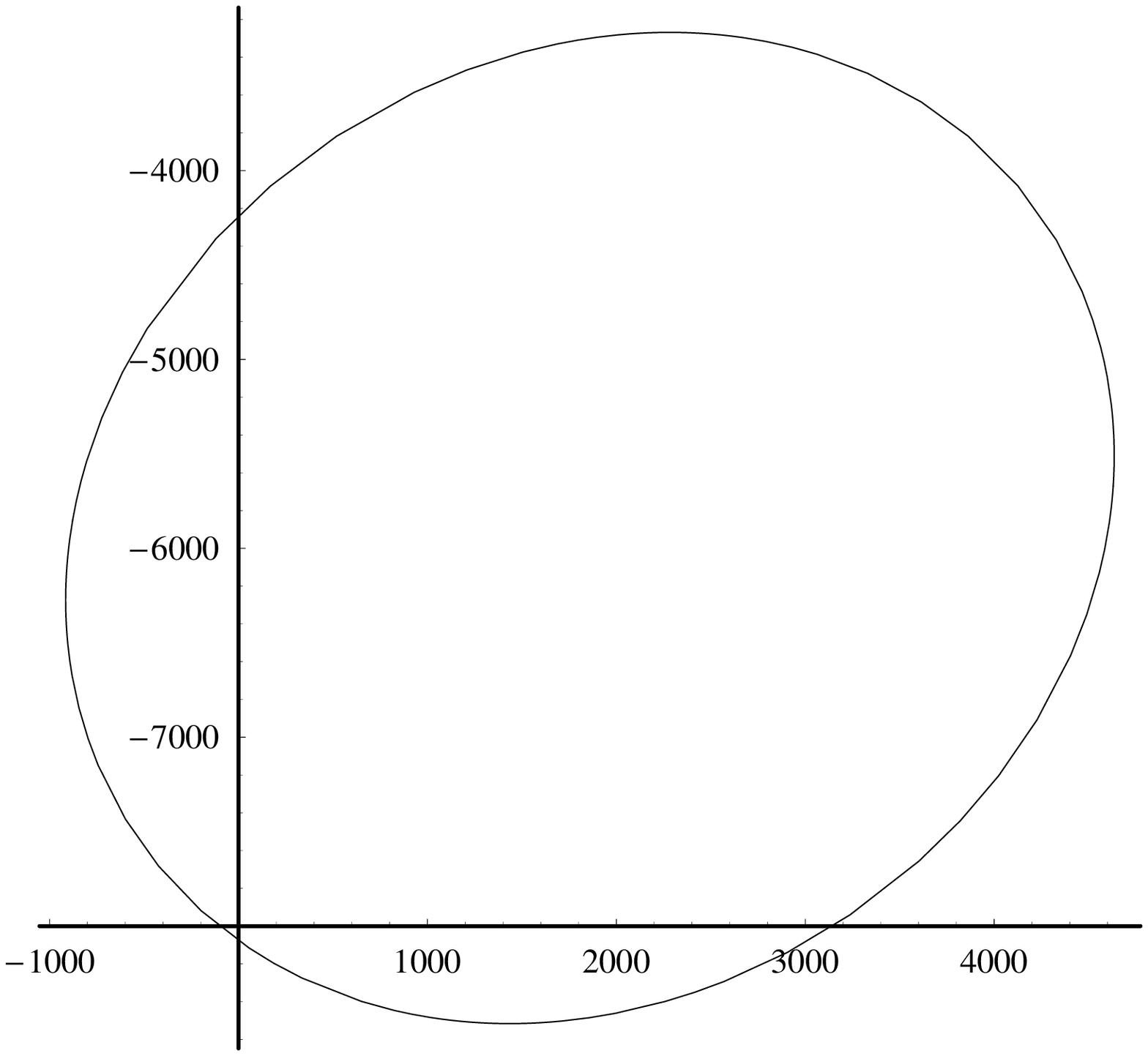}
\hspace{1cm}
\includegraphics[width=5cm]{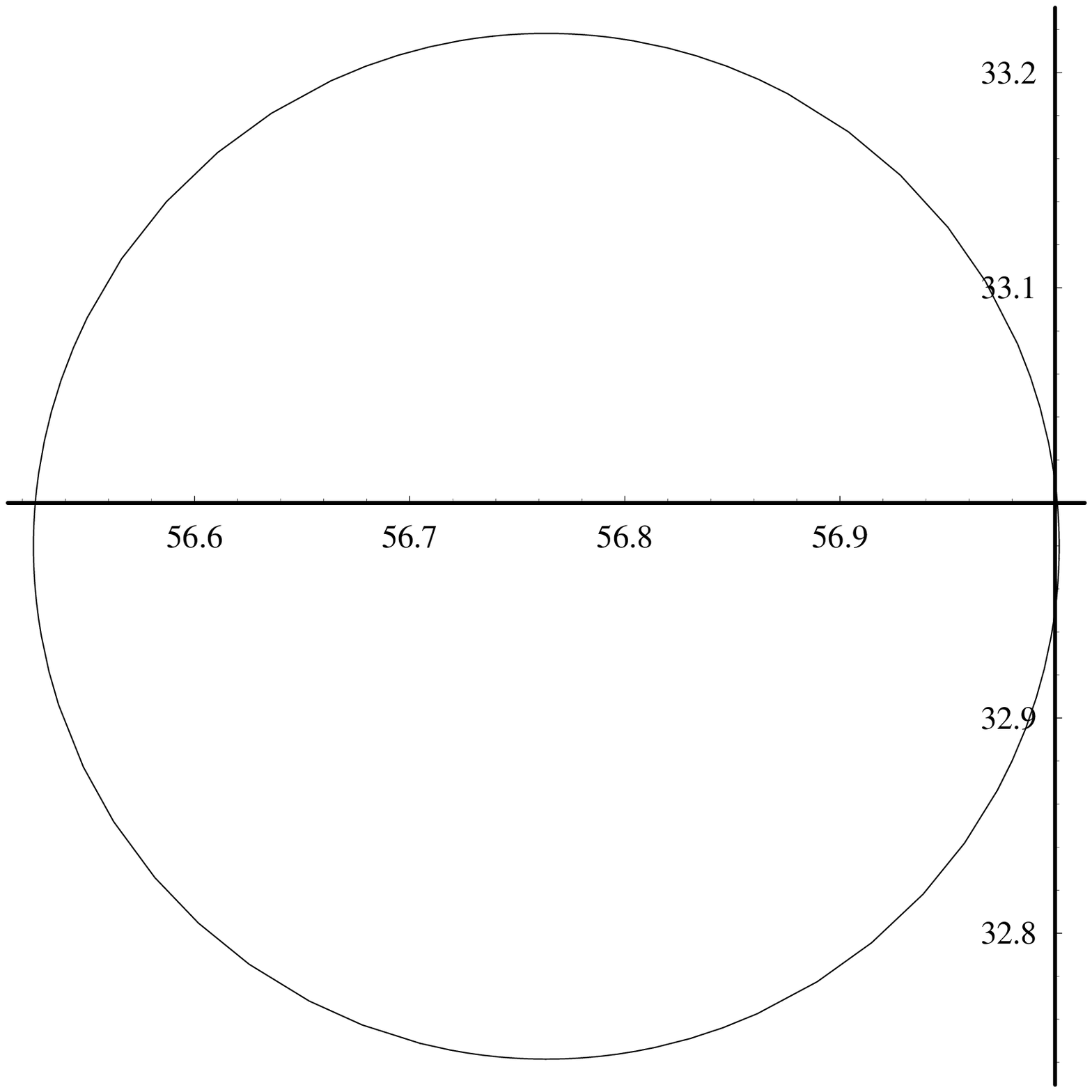}
\end{center}
\hspace{1.5cm}$\partial V(z_0,\lambda)$ \hspace{4.5cm}  $\partial V(z_0,\lambda)$
\caption{Region of variability for $\log f(z_0)$}
\end{figure}
\begin{center}
$\begin{array}{ll}
z_0 =0.737135+ 0.496542i       \hspace{3cm} & z_0 =-0.00588894 - 0.00496324i    \\
\lambda =-0.00646307 - 0.0167039i           & \lambda =-0.0472837 + 0.0970889i  \\
\mu =14038.5 + 9544.66i                     & \mu=25447.1- 2011.7i
\end{array}$
\end{center}

\begin{figure}[htp]
\begin{center}
\includegraphics[width=5cm]{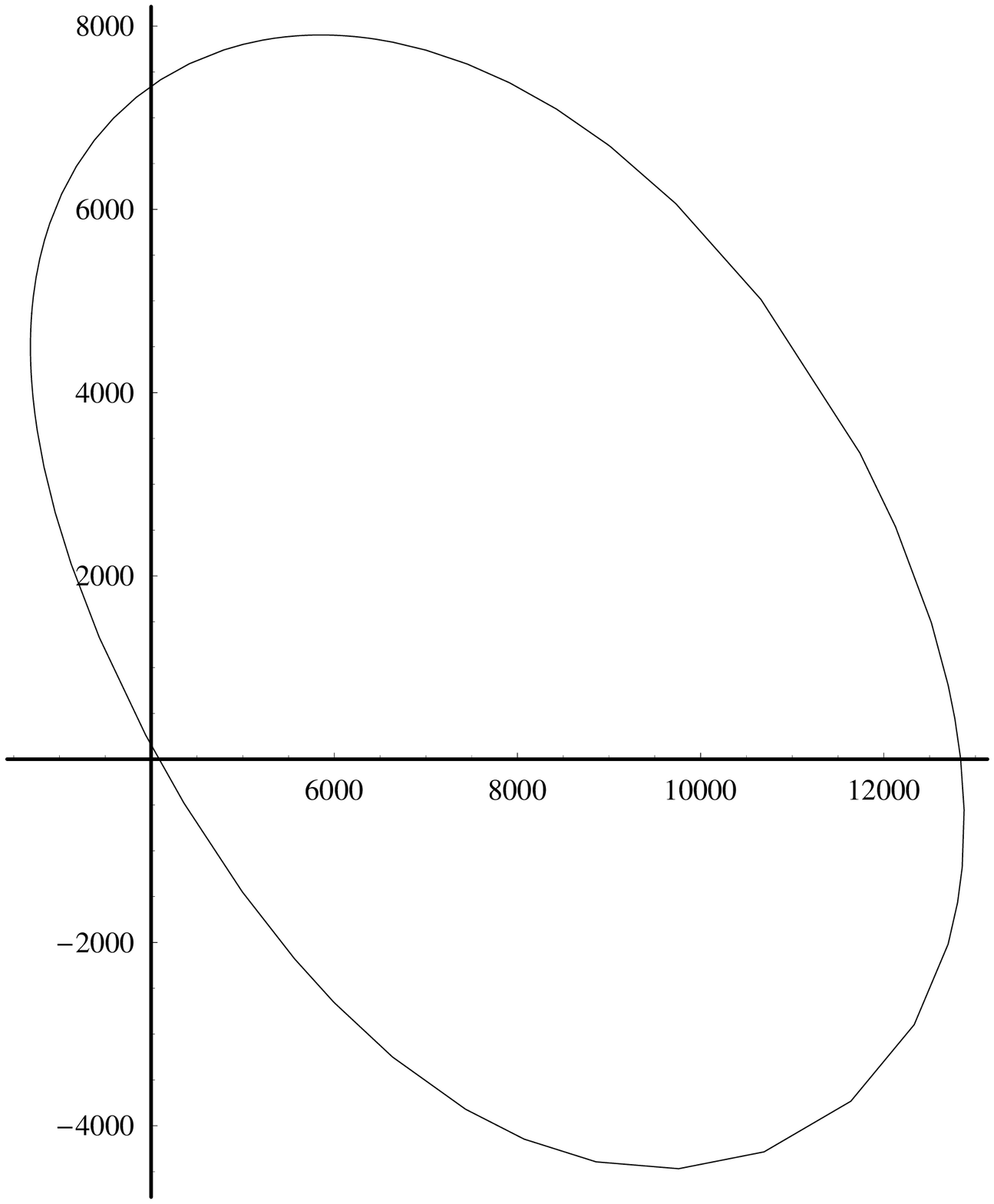}
\hspace{1cm}
\includegraphics[width=5cm]{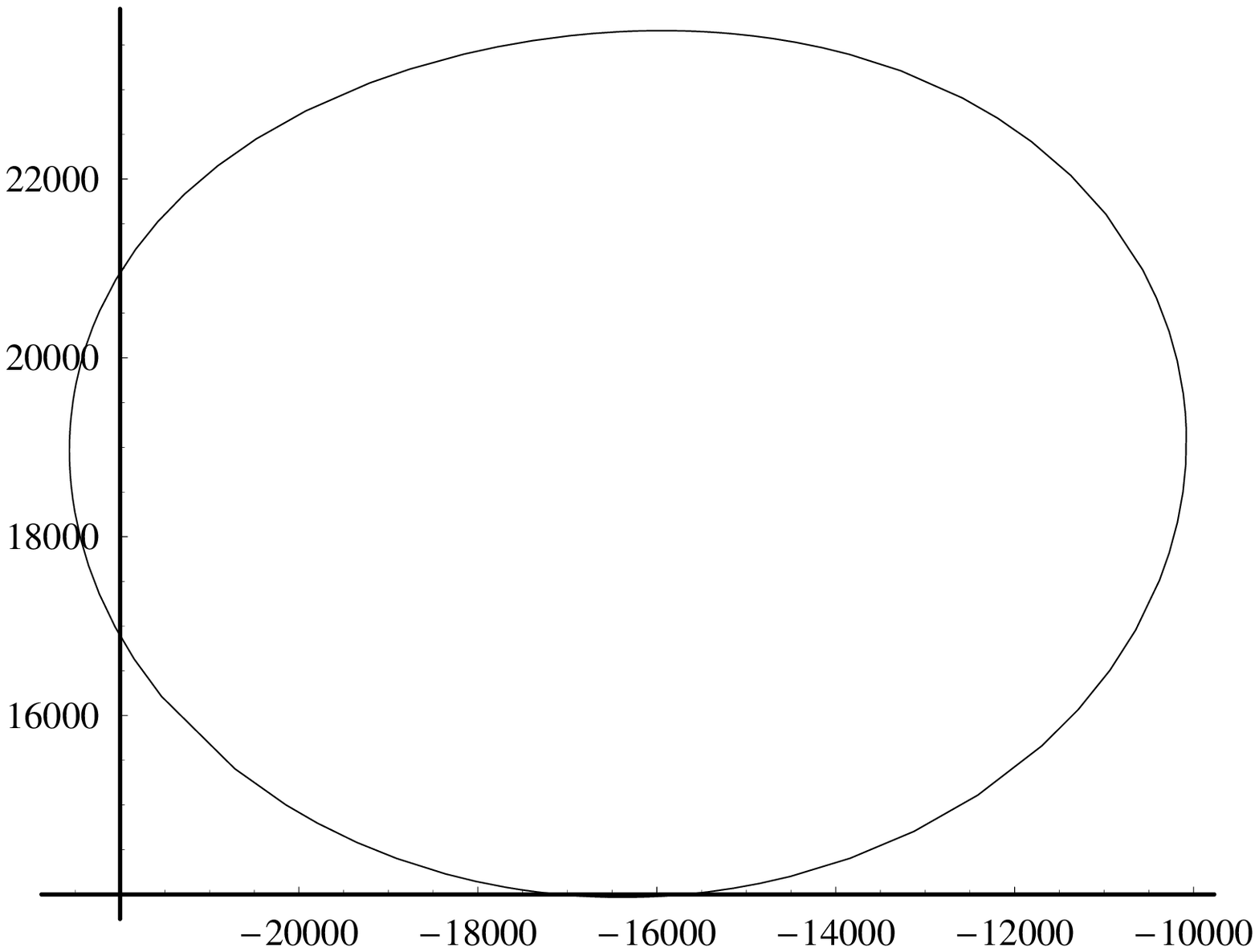}
\end{center}
\hspace{1.5cm}$\partial V(z_0,\lambda)$ \hspace{4cm}  $\partial V(z_0,\lambda)$
\caption{Region of variability for $\log f(z_0)$}
\end{figure}
\begin{center}
$\begin{array}{ll}
z_0 =0.556307- 0.814404i       \hspace{3cm} & z_0 =0.880992- 0.328223i    \\
\lambda =0.226895- 0.384635i                & \lambda =-0.0326596 + 0.656304i  \\
\mu =13589.3- 25797.8i                      & \mu= 39935.5+ 11412i
\end{array}$
\end{center}

\end{document}